\documentclass[11pt]{amsart}

\title{A simplified and generalized treatment of DES related ciphers}
\author{L. Babinkostova $^{1\P}$, A. M. Bowden$^{2*}$, A. M. Kimball$^{3*}$, and K. J. Williams$^{1*}$ }

\thanks{$^{\P}$ Communicating Author: liljanababinkostova@boisestate.edu}
\thanks{$^{*}$ Indicates Undergraduate Student}

\usepackage{amsfonts}
\usepackage{amsmath}
\usepackage{amsthm}
\usepackage{amssymb}
\usepackage{graphicx}
\usepackage{multirow}
\usepackage{setspace}
\usepackage{url}
\usepackage{subfigure}
\usepackage{color}

\newcommand{\naturals}{{\mathbb N}}
\newcommand{\integers}{{\mathbb Z}}

\newcommand{\KK}{\mathcal{K}}

\newcommand{\MM}{\mathcal{M}}
\newcommand{\TT}{\mathcal{T}}
\newcommand{\alt}{\mathcal{A}}
\newcommand{\sym}{\mathcal{S}}

\newcommand{\CC}{\mathcal{C}}

\newtheorem{definition}{{\bf Definition}}
\newtheorem{theorem}{{\bf Theorem }}
\newtheorem{lemma}[theorem]{{\bf Lemma }}
\newtheorem{corollary}[theorem]{{\bf Corollary}}

\newtheorem{conjecture}{{\bf Conjecture}}

\newcommand{\epf}{\diamondsuit} 

\graphicspath{%
    {converted_graphics/}
    {/}
}

\begin{document}
\maketitle
\begin{abstract}
This work is a study of DES-like ciphers where the bitwise exclusive-or (XOR) operation in the underlying Feistel network is replaced by an arbitrary group operation. We construct a two round simplified version of DES that contains all the DES components and show that its set of encryption permutations is not a group under functional composition, it is not a pure cipher and its set of encryption permutations does not generate the alternating group. We present a non-statistical proof that for $n\leq 6$ the set of $n$-round Feistel permutations over an arbitrary group do not constitute a group under functional composition.
  
\end{abstract}
\section{Introduction}

The Data Encryption Standard (DES), now replaced by AES (the Advanced Encryption Standard) as official symmetric key standard,
has been an official US government standard for more than 20 years. DES and its variants are still commonly used in electronic financial transactions, secure data communications, and the protection of passwords or PINs \cite{Mao}. Being the first commonly used modern block cipher, DES has undergone several thorough security analyses. Components of the DES architecture are still the fundamental building blocks for several contemporary symmetric key systems.

The full description of DES is given in \cite{NIST}, but we begin describing the details of DES that are relevant to this paper.
DES is a 16-round Feistel cipher acting on the space of 64-bit messages under the control of a 56-bit key. A description of the Feistel cipher is given in  \cite{F}. Each Feistel round is based on a round specific function $f$ mapping 32-bit strings to 32-bit strings. The function $f$ is derived from the round number and from an initial cryptographic key and involves both permutation and substitution operations. It has four components: 
\begin{enumerate}
  \item An expansion permutation $E$ which expands the right-half block $y$ from 32 bits to 48 bits.
  \item Key mixing which combines $E(y)$ with the subkey $k_i$ used in round $i$ by means of the bitwise-XOR operation.
  \item Substitution which divides the block $E(y)\oplus k_i$ into eight 6-bit blocks and each 6-bit block is transformed into  a 4-bit block according to a non-linear transformation, provided in the form of a lookup table (S-boxes).
  \item Permutation which rearranges the 32-bit output of the eight S-boxes according to a fixed permutation.
\end{enumerate}

Processing of a 64 bit block of data uses the bitwise exclusive-or (XOR) operation defined on 48-bit strings. The set of 48-bit strings, denoted $\{0,\,1\}^{48}$, endowed with the bitwise XOR operation denoted $\oplus$, is a group. This group is in fact the \emph{external direct product}\footnote{See Chapter 8 of \cite{G} for information on external direct products.} of 48 copies of the group $\integers_2$, which is the set $\{0,\,1\}$ endowed with the binary operation of modulo 2 addition. 

For each key $k$, the corresponding DES encryption operation is a function, say $T_k$, which maps 64-bit data blocks to 64-bit blocks, and thus $T_k$ is a function from $\{0,\,1\}^{64}$ to $\{0,\,1\}^{64}$. For each key $k$, $T_k$ is a one-to-one (and thus onto) function from $\{0,\,1\}^{64}$ to $\{0,\,1\}^{64}$, and thus a permutation of the set $\{0,\,1\}^{64}$. We shall call the DES encryption functions \emph{DES permutations}. The set of DES permutations is a subset of the symmetric group ${\sf S}_{2^{64}}$ which has $(2^{64})!$ elements. 

If the set of DES permutations were closed under functional composition (and thus a subgroup of ${\sf S}_{2^{64}}$) then multiple encryption using several DES keys would be equivalent to a single encryption by a single DES key. In this case security features of systems like Triple-DES would not exceed that of DES. It has been known since the 1990's that the set of DES encryption permutations does not constitute a group under the operation of functional composition. In \cite{CW} the authors presented a statistical test to show that the indexed set of permutations is not closed under functional composition and so is not a group. 

Our work is a study of cryptosystems with DES-like architecture, but for which the bitwise exclusive-or (XOR) operation in the underlying Feistel network is replaced by a binary operation in an arbitrary finite group. The idea of replacing the XOR operation in DES with another operation has been considered by several authors: Biham and Shamir \cite{BS} show that replacing some of the XOR operations in DES with additions modulo $2^n$, makes their differential attack less powerful. Carter, Dawson, and Nielsen \cite{CDN} show a similar phenomenon when the XOR operation in DES is replaced by addition using a particular Latin Square. Patel, Ramzan and Sundaram \cite {XOR} studied  Luby-Rackoff ciphers over arbitrary finite groups. They constructed a four round Luby-Rackoff cipher, operating over finite groups of characteristic greater than 2, and showed that such a cipher is secure against adaptive chosen plaintext and adaptive chosen ciphertext attacks, has better time/space complexity and uses fewer random bits than the previously considered Luby-Rackoff ciphers based on the group $\integers_2$ with an XOR operation. As with DES, the Luby-Rackoff cipher involves the use of Feistel permutations independently keyed with pseudorandom functions. 

Before our 
investigation there was no information on whether the set of encryption permutations of such systems based on a different finite group operation is closed under functional composition. We show 
that, in non-pathological cases, for $n\leq 6$ the set of encryption permutations generated by $n$-round Feistel permutations is not closed under functional composition. We must also note that we found a fairly simple deductive proof of this fact instead using statistical tests as in \cite{CW}.

Knowing the order of the group generated by the encryption permutations is also an important algebraic question about the security of the cryptosystem. Coppersmith and Grossman have shown \cite{CE} that in principle DES-like components can generate any permutation of the alternating group $\mathcal{A}_{2^{64}}$ (all even permutations, i.e. those that can be represented by an even number of transpositions). In 1983 S. Even and O. Goldreich showed that DES-like functions are contained within the alternating group \cite {EG}. Furthermore, in 1998 R. Wernsdorf \cite {W} showed that the one-round encryption permutations of DES generate the alternating group, $\mathcal{A}_{2^{64}}$. It is still not known whether 16-round DES permutations generate the alternating group. Using the special properties of the so called weak keys it has been shown in \cite {CO} that the set of DES permutations generates a very large group, with a lower-bound of $2^{2499}$ for its size.

In \cite {BBKW} it is 
shown how the replacement of the XOR operation in the underlying Feistel network of DES with binary operation in an arbitrary finite group can affect which group is generated by the $n$-round DES permutations. 

Since any direct analysis of DES is computationally intensive it is sometimes not feasible to directly analyze DES. Thus, simplified analogs of DES have been introduced. Such simplified versions of DES were introduced by E. Schaeffer (called S-DES) in \cite{SDES}, and W. Trappe and L. Washington (called B-DES) in \cite{BDES}. Both versions simulate the basic architecture of DES. As with DES, the fundamental computational structure underlying these simplified versions of DES is the group $\integers_2$ endowed with the XOR operation. In \cite{KT}, J. Konikoff and S. Toplosky, showed that the group of permutations generated by these simplified versions of DES is in the case of S-DES the alternating group on 256 elements and in the case of B-DES the alternating group on 4096. As with DES, there was no information on whether or how the security of corresponding analogs of S-DES or B-DES is affected when the XOR operation is replaced with another group operation.

In this paper we introduce a new simplified version of DES, which we call E-DES. The main innovation in E-DES is that we base it on the finite group ${\mathbb Z}_3$ with the operation of modulo 3 addition. We show that its set of encryption permutations does not form a group under functional composition. We show this both by giving a deductive proof of the statement, and by using the cycling closure test as was used in the case of DES for the same question. Also, we show that the group generated by the E-DES permutations is not the alternating group as it contains odd permutations. We gave careful attention to the known original design criteria for components of DES, particularly the substitution S-boxes, when we designed E-DES. We also show that this cryptosystem is not pure, another important property for the security of any cryptosystem. 

Our construction of E-DES based on $\integers_3$ can be easily adjusted to develop a cipher based on a finite group different from ${\mathbb Z}_3$ while preserving security features of E-DES. In particular, one can use Elliptic Curve groups in anticipation that several computationally hard problems for these groups may be used to further enhance the security of a DES-like cryptosystem based on these groups. An attempt for use of the Elliptic Curve groups in DES was made in \cite {AHE}, but unfortunately not as replacement of the bitwise operation in the underlying Feistel network.  

This paper is organized as follows: In Section 2 we state several definitions and give a general description of block ciphers based on multiple round Feistel networks. In Section 3 we describe a technique of constructing Feistel functions. In Section 4 we give the specifics of the design of our instance of E-DES and give an example of an encryption in this system. In Section 5 we describe the cycling closure test and how it can be used to address various questions about the algebraic structure of any finite cryptosystem. Our results based on these tests address the question of whether the E-DES encryption permutations  constitute a group under functional composition, whether the cryptosystem is pure and which group is generated by such cryptosystems. In Section 6 we address these and other group theoretic properties concerning of Feistel based block ciphers in general. 

\section {Definitions and Notation}

We follow the notation and terminology of \cite {Kaliski}. A cryptosystem is an ordered 4-tuple $(\MM,\,\CC,\,\KK,\,T)$ where $\MM$, $\CC$, and $\KK$ are called the {\it message space}, the {\it ciphertext space}, and the {\it key space} respectively, and where $T: \MM\times\KK\rightarrow \CC$ is a transformation such that for each $k\in\KK$, the mapping $T_k:\MM\rightarrow\CC$ is invertible.

For any cryptosystem $\Pi=(\MM,\,\CC,\,\KK,\,T)$, let $\TT_{\Pi}=\{T_k:k\in\KK\}$ be the set of all encryption transformations. In addition, for any transformation $T_k\in\TT$, let $T_k^{-1}$ denote the inverse of $T_k$.  In a cryptosystem where $\MM=\CC$ the mapping $T_k$ is a permutation of $\MM$. We consider only cryptosystems for which $\MM = \CC$. The set of all permutations of the set $\MM$ is denoted $\sym_{\MM}$. Under the operation of functional composition $\sym_{\MM}$ forms a group called \emph{the symmetric group} over $\MM$. The symbol $\langle\TT_{\Pi}\rangle$ denotes the subgroup of $\sym_{\MM}$ that is generated by the set $\TT_{\Pi}$. 

A cryptosystem $\Pi$ is called {\it closed} if its set $\TT_{\Pi}$ of encryption transformations is closed under functional composition i.e for every $k_1, k_2\in \KK$ there is $k_3\in\KK$ such that $T_{k_1}T_{k_2}=T_{k_3}$. By Theorem 3.3 from \cite{G} $\Pi$ is closed if and only if its set of encryption transformations $\TT_{\Pi}$ is a group under functional composition.  Thus, the cryptosystem $\Pi$ is closed if, and only if, $\TT_{\Pi} = \langle \TT_{\Pi}\rangle$. 

In \cite{S} Shannon generalized the idea of closed cipher. A cryptosystem is {\it pure} if and only if for every three keys $k_1, k_2,$ and $k_3$ there exists a key $k_4$ such that $T_{k_1}T_{k_2}^{-1}T_{k_3} = T_{k_4}$. 
One can show that $\Pi$ is pure if and only if for every $T\in\TT_{\Pi}$ the set $T^{-1}\TT_{\Pi} = \{T^{-1}T_k:k\in\KK\}$ forms a group under functional composition. It is known that every closed cryptosystem is pure, but not every pure cryptosystem is closed.  

To analyze the algebraic properties of Feistel based ciphers it is also useful to introduce the following definitions about permutation groups. For any subgroup $S\subseteq \sym_{\MM}$, for any $m\in\MM$, the set $orb_S(m)=\{\phi(m):\phi\in S\}$ is called the \emph{orbit} of $m$ under $S$. The set $stab_S(m)=\{\phi\in S:\phi (m)=m\}$ is called the \emph{stabilizer} of $m$ in $S$. In Section 5 we will make use of the following well-known theorem.
\begin{theorem}
Let $S$ be a finite group of permutations of a set $M$. Then for any $m\in M$, 
\[
\mid S\mid=\mid orb_S(m)\mid \cdot \mid stab_S(m)\mid
\]
\end{theorem}

\vspace{0.05in}
\section {Feistel Networks}

A multi-round block cipher is a cipher involving the sequential application of similar invertible transformations (called round functions or round transformations) to the plaintext. All round transformations are usually key-dependent and the transformation of round $i$ obtains its own subkey $k_i$ which is derived from the cipher key $k$ using a key-schedule algorithm. {\it Feistel networks} constitute an important design principle underlying many block ciphers, including DES. They were described first by Horst Feistel during his work at IBM on the cipher Lucifer \cite {F}.  
 
\begin{definition}
Let $(G,\oplus)$ be a finite group. For a function $f:G^{t} \rightarrow G^{t}$ the function $\sigma_f: G^{2t} \rightarrow G^{2t}$defined by $\sigma_f(x,y) = (y, x \oplus f(y))$ is called  a {\it Feistel function}.
\end{definition}

A {\it Feistel network} consists of repeated applications of Feistel functions with different round functions $f$ used in each round. By definition, a Feistel network with $n$-rounds is the permutation function 
\[
\Psi^{n}_{2t}(f_1,f_2,\cdots,f_n)= \sigma_{f_1}\circ\sigma_{f_2}\circ\cdots\circ\sigma_{f_n}
\]

Typically, the round functions are chosen to be highly nonlinear key-dependent functions with good diffusion. A standard way to provide these properties for the round functions is to use a substitution structure (S-boxes).

Let $F_t(G)$ be the set of all functions $f:G^t\rightarrow G^t$. The set $F_t(G)$ admits some operations. For $f$ and $g$ in ${\sf F}_t(G)$ we define $f\odot g$ so that for each $y\in G^t$,
\[
  (f\odot g)(y) = f(y) \oplus  g(y).
\]
The symbol ``$\oplus$" above denotes the group operation of the product group $G^t$. Now $({\sf F}_t(G),\, \odot)$ is a group. Additionally, the functional composition operation $\circ$ on ${\sf F}_t(G)$ will be featured.

A {\it random Feistel network} with $n$ rounds, is a Feistel network in which the round functions $f_1,\cdots, f_n$ are  randomly and independently chosen functions form the set $F_t(G)$. These networks are also known as ``Luby-Rackoff constructions with $n$ rounds". 

\begin{lemma}
Let $(G,\oplus)$ be a finite group. For each function $f\in {\sf F}_t(G)$ the Feistel function $\sigma_f$ is a permutation of the set $G^t\times G^t$.
\end{lemma}

Thus, for each function $f\in {\sf F}_t(G)$, the Feistel function $\sigma_f$ is a member of $\sym_{\mid G^{2t}\mid}$, the permutation group of the finite set $G^t\times G^t$.

\subsection{Feistel functions derived from S-boxes}

In several practical implementations of Feistel network based cryptosystems, including classical DES, the members of $F_t(G)$ are constructed in a very specific way from an input key parameter, and a selected set of substitution tables called {\it S-boxes}. Much of the security of a block cipher based on a Feistel network depends on the properties of the substitution boxes (S-boxes) used in the round function.  
The DES S-boxes are reported to have been designed to conform to a number of criteria as they are the part of the system where the cipher function gets its security. For more details on the properties of the S-boxes used in DES, see for example \cite {CO}, \cite{BGK}, \cite {BMP} and  \cite {CCG}. In this section we describe the design of functions $f:G^{t} \rightarrow G^{t}$ where $G$ is some finite group, from S-boxes. 

An S-box is a lookup table with $k=\vert G\vert^i$ rows and $m=\vert G\vert^j$ columns. The entries of an S-box will be conceived of as $j$-nit sequences over $G$. Let $n$ such S-boxes be given. The function $f$ constructed from these $n$ S-boxes must be a function from $G^t$ to $G^t$. Thus, the inputs of the function $f$ will be a $t$-nit \footnote{here, nit (from ``$n$-ary digit'', analogous to bit) means group element} sequence of group elements. 

This input is used to construct an output from the $n$ S-boxes by reading off from this input an S-box number and row-number and column-number for that S-box, and then using the entries in these positions of the indicated S-boxes to construct the output. There are several approaches to indicating the S-boxes. 

We require that all $n$ S-boxes are ``active" in constructing the value of the function $f$. Thus, from the $t$-nit input we must read for each of the $n$ S-boxes the corresponding row- and column-information. Note that from our specification of the S-box dimensions above, the row number can be coded as an $i$-nit string of elements of $G$, while the column number can be coded as a $j$-nit string of elements of $G$. Thus, to specify a row and column we use an $(i+j)$-nit string, plus a convention indicating which $i$ of these nits to use in which specific ordering to select the row, and what ordering on the remaining $j$-nits to use to select the column. Let such a row/column number convention be fixed. 
To keep the selections from the $S$-boxes independent of each other, we use an $(i+j)\hspace{-0.04in}\cdot \hspace{-0.04in}n$-nit string derived from the input element of $G^t$, a $t$-nit string. The output is obtained by concatenating the $n$ input-indicated $j$-nit entries from each of the S-boxes, thus obtaining a $(j\hspace{-0.04in}\cdot\hspace{-0.04in} n)$-nit output. Thus, we require that $t=j\hspace{-0.04in}\cdot\hspace{-0.04in} n$. 

Since $t <(i+j)\hspace{-0.04in}\cdot \hspace{-0.03in}n$, we expand the $t$-nit input using a carefully chosen expansion function 
\[
  E: G^t \longrightarrow G^{(i+j)\hspace{-0.01in}\cdot \hspace{-0.01in}n}.
\]

The other input parameter to the S-box is an $(i+j)\hspace{-0.04in}\cdot \hspace{-0.04in}n$-nit key, $K$. This key and $E(R)$, the expansion of the $t$-nit string $R$, is then used to construct the input to the S-box. The input to the S-box is $K \oplus E(R)$ where the operation $\oplus$ is the nit-wise group operation on the product group $G^{(i+j)\hspace{-0.03in}\cdot \hspace{-0.03in}}$.
The $(i+j)\hspace{-0.04in}\cdot \hspace{-0.04in}n$-nit quantity $ K \oplus E(R)$ is separated into $n$ blocks of consecutive $(i+j)$-nits each, with string number $s$ designated as the row/column selection code for S-box number $s$. 

The output is obtained by concatenating the $j$-nit entries from each of the S-boxes in canonical order to obtain a $t$-nit ($t=j\hspace{-0.03in}\cdot\hspace{-0.03in} n)$ output. Figure \ref{S-boxPicture} illustrates this construction based on three S-boxes.

\begin{figure}[tbp] 
  \centering  
   \includegraphics[width=2.7in,height=2.54in,keepaspectratio]{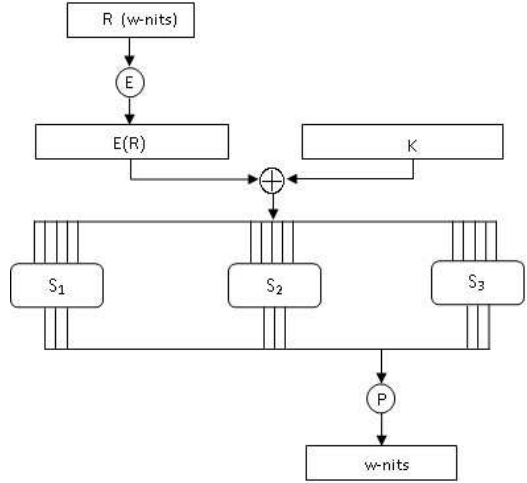}
  \caption{Calculation of $f(R,K)$}
  \label{Calculation of $f(R,K)$}
\end{figure}
\section {A simplified DES cipher based on $\integers_3$}

\begin{figure}[tbp] 
  \centering
    \includegraphics[bb=0 0 262 321,width=2.25in,height=2.76in,keepaspectratio]{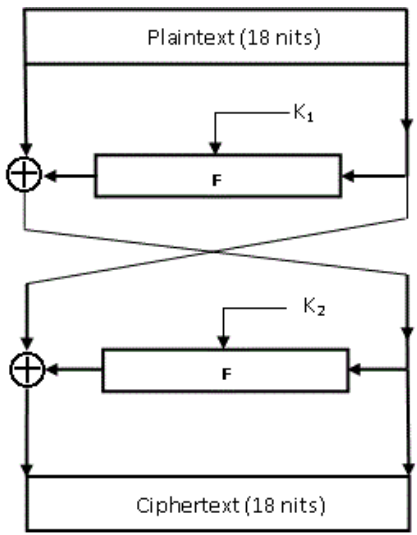}
    \caption{E-DES structure}
  \label{E-DES structure}
\end{figure}

We define a simplified version of DES, called E-DES, by declaring the message-space $\MM$ and the ciphertext space $\CC$ to be $\MM=\CC=\{0,1,2\}^{18}$, and by declaring the key space $\KK$ to be $\KK=\{0,1,2\}^{20}$.  Figure \ref{E-DES structure} illustrates the overall structure of this simplified DES cipher and Figure \ref{S-boxPicture} illustrates the round function $f$.

The inputs for the E-DES encryption algorithm are an $18$-nit \footnote{here, nit (from ``$n$-ary digit'', analogous to bit) means group element} block of plaintext (example: 110120121220120121) and a 20-nit key. The algorithm produces an 18-nit block of ciphertext as output. For a fixed key $k$ the function $f_k$ derived from the S-boxes as described before takes as input the data passing through the encryption algorithm and a 15-nit subkey. The mapping $T_k$ defined from the key $k$ by the E-DES encryption algorithm involves the sequential application of functions and can be expressed as 
\[
   T_k=P^{-1}\circ \theta\circ \sigma_{f_{k_2}}\circ\sigma_{f_{k_1}}\circ P. 
\]
where $\theta$ is the ``swap" function 
           \begin{center}
               $\theta(x,y)=(y,x)$
            \end{center}
from $\integers_3^{18}$ to $\integers_3^{18}$. 

Note that $\Theta = \Theta^{-1}$. Recall that $(f\circ g)^{-1}= g^{-1}\circ f^{-1}$. Then the decryption algorithm can be expressed as
\[
   T^{-1}_k=P^{-1}\circ \sigma^{-1}_{f_{k_1}}\circ\sigma^{-1}_{f_{k_2}}\circ\theta\circ P. 
\]
The inverse of $\sigma^{-1}_{f_{k_i}}$, $i=1, 2$ is   
\[
\sigma^{-1}_{f_{k_i}}(x,y) = (x \ominus f(y), y),
\]
where $\ominus$ is the operation subtraction modulo 3, the inverse group operation in $\integers_3$.
 
The quantity $k_1$ is a 15-nit subkey of the input key $k$ and is derived from $k$ using the compression permutation 
\[
  CP_1={
       \left(\begin{smallmatrix}
               16 & 17 & 12 & 15 & 20 & 10 & 11 & 3 & 7 & 19 & 13 & 9 & 8 & 1 & 18\\
             \end{smallmatrix}
       \right)}
\]
Similarly, $k_2$ is a subkey of $k$ obtained from $k$ by using the compression permutation 
\[
  CP_2={
       \left(\begin{smallmatrix}
               6 & 7 & 2 & 20 & 4 & 3 & 9 & 8 & 18 & 10 & 15 & 14 & 11 & 12 & 5
             \end{smallmatrix}
       \right)}
\]

We shall now describe each of the functions constituting the encryption permutation $T_k$. Following this description, we illustrate the algorithm by going through the individual steps using an explicit input plaintext and an explicit input key.

\subsection{Initial and Final Permutations}
The input to the algorithm is an 18-nit block of plaintext which is first permuted using the permutation

\[ 
  P={
     \left(
            \begin{smallmatrix}
                1 &  2 & 3  &  4 & 5 &  6 &  7 & 8 & 9  & 10 & 11 & 12 & 13 & 14 & 15 & 16 & 17 & 18\\
               11 & 12 & 2  & 13 & 9 &  1 &  5 & 8 & 16 & 17 &  4 & 18 & 15 &  7 & 10 &  3 &  6 & 14
            \end{smallmatrix}
     \right).}
\]

For example, the initial permutation moves nit in position 6 of the plaintext to position 1, the nit in position 3 of the plaintext to position 2, the nit in position 16 of the plaintext to position 3, and so forth.

The final permutation is the inverse of $P$, and is the permutation
\[
  P^{-1}={
     \left(
            \begin{smallmatrix}
                1 & 2 & 3  &  4 & 5 &  6 &  7 & 8 & 9 & 10 & 11 & 12 & 13 & 14 & 15 & 16 & 17& 18\\
                6 & 3 & 16 & 11 & 7 & 17 & 14 & 8 & 5 & 15 & 1  & 2  & 4  & 18 & 13 & 9 & 10 & 12
            \end{smallmatrix}
     \right).}
\]

The initial permutation $P$ and the corresponding final permutation $P^{-1}$ do not affect the security of the cryptosystem. 

\subsection{The Expansion Permutation}
This operation expands the right half of the 18-nit data-block being processed from 9-nits to 15-nits by changing the order of the nits as well as repeating certain nits. This operation has two properties: It makes the right half the same size as the subkey for the $\oplus_{\,mod\,\,3}$ operation  and it provides the appropriate length nit-sequence for use during the substitution operation. The expansion permutation is given by 
\[
  E={
      \left(
         \begin{smallmatrix}
             9 & 1 & 2 & 3 & 4 & 5 & 6 & 3 & 4 & 5 & 6 & 7 & 8 & 9 & 1
         \end{smallmatrix}
      \right).}
\]   

\subsection{The S-box substitution}
The most fundamental encryption step in E-DES, directly impacting the security of E-DES, is the application of the substitution boxes, or S-boxes. There are three different S-boxes. Table 2 shows these three S-boxes. To achieve the ``non-linearity" property  and the randomness of the output we used the following criteria in the design of the S-boxes:
\begin{enumerate} 
  \item No S-box is a linear or affine function of the input.  
  \item Each ``row" of an S-box contains all possible outputs.
\end{enumerate}
Note that the entries displayed in Table 2 are ordinary integers between 0 and 26. Each S-box transforms a 5-nit input to a 3-nit output as follows: The 15-nit result of the expansion permutation is divided into three 5-nit sub-blocks. Each separate block is the input for a separate S-box: The first block is input for S-box 1, the second block is input for S-box 2 and the third block is input for S-box 3. 

Each S-box is a table of 9 rows and 27 columns. The rows are numbered by 0 through 8, while the columns are numbered by 0 through 26. Each entry in the 9-by-27 S-box is, when expressed in base 3, a 3-nit number. The 5-nit input of an S-box specifies the row and column number of the S-box entry that is the output for that input. This is done as done as follows: Let $n_1, n_2, n_3, n_4$, and $n_5$ be the 5 nits of input listed in order of occurrence in the input. Nits $n_1$ and $n_5$ are combined to form 2-nit number in base 3, corresponding to one of the decimal numbers from 0 to 8: This decimal number specifies a row number in the S-box under consideration. The middle 3 nits, $n_2, n_3, n_4$ are combined to form a 3-nit number in base 3, corresponding to one of the decimal numbers from 0 to 26: This decimal number specifies a column in the S-box under consideration. Here 00 corresponds to row 1, and 000 corresponds to column 1. 

For example, suppose that the input to the second S-box is 22010. The first and the last nit combine to form 20, which corresponds to row number 6, which is by our convention the seventh row, of the second S-box. The middle 3 nits combine to form 201, which correspond to the number 19, indicating by our convention the 20-th column of the same S-box. The entry at the intersection of the seventh row and twentieth column of S-box 2 is 11. Since 102 is the base 3 representation of 11, the 3-nit value 102 is the output from S-box 2, given the 5-nit input 22010. 

\begin{figure}
\begin{center}
{ \Large Table 1: {\bf S-boxes}}
\line(1,0){350}
\end{center}
{\bf S-box 1}
\[
\begin{smallmatrix}
24 & 25 & 6 & 16 & 3 & 7 & 1 & 18 & 26 & 5 & 10 & 9 & 19 & 23 & 13 & 12 &15 & 8 & 20 & 17 & 2 & 11 & 0 & 21 & 14 & 4 & 22\\
&&&&&&&&&&&&&&&&&&&&&&&&&&&\\
17 & 18 & 26 & 9 & 23 & 0 & 21 & 11 & 19 & 25 & 3 & 2 & 12 & 16 & 6 & 5 & 8 & 1 & 13 & 10 & 22 & 4 & 20 & 14 & 7 & 24 & 15\\
&&&&&&&&&&&&&&&&&&&&&&&&&&&\\
16 & 17 & 25 & 8 & 22 & 26 & 20 & 10 & 18 & 24 & 2 & 1 & 11 & 15 & 5 & 4 & 7 & 0 & 12 & 9 & 21 & 3 & 19 & 13 & 6 & 23 & 14\\
&&&&&&&&&&&&&&&&&&&&&&&&&&&\\
10 & 11 & 19 & 2 & 16 & 20 & 14 & 4 & 12 & 18 & 23 & 22 & 5 & 9 & 26 & 25 & 1 & 21 & 6 & 3 & 15 & 24 & 13 & 7 & 0 & 17 & 8\\
&&&&&&&&&&&&&&&&&&&&&&&&&&&\\
21 & 22 & 3 & 23 & 0 & 4 & 25 & 15 & 23 & 2 & 7 & 6 & 16 & 20 & 10 & 9 & 12 & 5 & 17 & 14 & 26 & 8 & 24 & 18 & 11 & 1 & 19\\
&&&&&&&&&&&&&&&&&&&&&&&&&&&\\
26 & 0 & 8 & 18 & 5 & 9 & 3 & 20 & 1 & 7 & 12 & 11 & 21 & 25 & 15 & 14 & 17 & 10 & 22 & 19 & 4 & 13 & 2 & 23 & 16 & 6 & 24\\
&&&&&&&&&&&&&&&&&&&&&&&&&&&\\
3 & 4 & 12 & 22 & 9 & 13 & 7 & 24 & 5 & 11 & 16 & 15 & 25 & 2 & 19 & 18 & 21 & 14 & 26 & 23 & 8 & 17 & 6 & 0 & 20 & 10 & 1\\
&&&&&&&&&&&&&&&&&&&&&&&&&&&\\
5 & 6 & 14 & 24 & 11 & 15 & 9 & 26 & 7 & 13 & 18 & 17 & 0 & 4 & 21 & 20 & 23 & 16 & 1 & 25 & 10 & 19 & 8 & 2 & 22 & 12 & 3\\
&&&&&&&&&&&&&&&&&&&&&&&&&&&\\
11 & 12 & 20 & 3 & 17 & 21 & 15 & 5 & 13 & 19 & 24 & 23 & 6 & 10 & 0 & 26 & 2 & 22 & 7 & 4 & 16 & 25 & 14 & 8 & 1 & 18 & 9\\
\end{smallmatrix}
\]\vspace{.01in}

{\bf S-box 2}
\[
\begin{smallmatrix}
1 & 2 & 10 & 20 & 7 & 11 & 5 & 22 & 3 & 9 & 14 & 13 & 23 & 0 & 17 & 16 & 19 & 12 & 24 & 21 & 6 & 15 & 4 & 25 & 18 & 8 & 26\\
&&&&&&&&&&&&&&&&&&&&&&&&&&&\\
25 & 26 & 7 & 17 & 22 & 8 & 2 & 19 & 0 & 6 & 11 & 10 & 20 & 24 & 14 & 13 & 16 & 9 & 21 & 18 & 3 & 12 & 1 & 22 & 15 & 5 & 23\\
&&&&&&&&&&&&&&&&&&&&&&&&&&&\\
14 &  15 & 23 & 6 & 20 & 24 & 18 & 8 & 16 & 22 & 0 & 26 & 9 & 13 & 3 & 2 & 5 & 25 & 10 & 7 & 19 & 1 & 17 & 11 & 4 & 21 & 12\\
&&&&&&&&&&&&&&&&&&&&&&&&&&&\\
9 & 10 & 18 & 1 & 15 & 19 & 13 & 3 & 11 & 17 & 22 & 21 & 4 & 8 & 25 & 24 & 0 & 20 & 5 & 2 & 14 & 23 & 12 & 6 & 26 & 16 & 7\\
&&&&&&&&&&&&&&&&&&&&&&&&&&&\\
23 & 24 & 5 & 15 & 2 & 6 & 0 & 17 & 25 & 4 & 9 & 8 & 18 & 22 & 12 & 11 & 14 & 7 & 19 & 16 & 1 & 10 & 26 & 20 & 13 & 3 & 21\\
&&&&&&&&&&&&&&&&&&&&&&&&&&&\\
2 & 3 & 11 & 21 & 8 & 12 & 6 & 23 & 4 & 10 & 15 & 14 & 24 & 1 & 18 & 17 & 20 & 13 & 25 & 22 & 7 & 16 & 5 & 26 & 19 & 9 & 0\\
&&&&&&&&&&&&&&&&&&&&&&&&&&&\\
18 & 19 & 0 & 10 & 24 & 1 & 22 & 12 & 20 & 26 & 4 & 3 & 13 & 17 & 7 & 6 & 9 & 2 & 14 & 11 & 23 & 5 & 21 & 15 & 8 & 25 & 16\\
&&&&&&&&&&&&&&&&&&&&&&&&&&&\\
15 & 16 & 24 & 7 & 21 & 25 & 19 & 9 & 17 & 23 & 1 & 0 & 10 & 14 & 4 & 3 & 6 & 9 & 13 & 8 & 20 & 2 & 18 & 12 & 5 & 22 & 13\\
&&&&&&&&&&&&&&&&&&&&&&&&&&&\\
2 & 23 & 4 & 14 & 1 & 5 & 26 & 16 & 24 & 3 & 8 & 7 & 17 & 21 & 11 & 10 & 13 & 16 & 18 & 15 & 0 & 9 & 25 & 19 & 12 & 2 & 20\\
\end{smallmatrix}
\]\vspace{.01in}

{\bf S-box 3}
\[
\begin{smallmatrix}
4 & 5 & 13 & 23 & 10 & 14 & 8 & 25 & 6 & 12 & 17 & 16 & 26 & 3 & 20 & 19 & 22 & 15 & 0 & 24 & 9 & 18 & 7 & 1 & 21 & 11 & 2\\
&&&&&&&&&&&&&&&&&&&&&&&&&&&\\
6 & 7 & 15 & 25 & 12 & 16 & 10 & 0 & 8 & 14 & 19 & 18 & 1 & 5 & 22 & 21 & 24 & 17 & 2 & 26 & 11 & 20 & 9 & 3 & 23 & 13 & 4\\
&&&&&&&&&&&&&&&&&&&&&&&&&&&\\
7 & 8 & 16 & 26 & 13 & 17 & 11 & 1 & 9 & 15 & 20 & 19 & 2 & 6 & 23 & 22 & 25 & 18 & 3 & 0 & 12 & 21 & 10 & 4 & 24 & 14 & 5\\
&&&&&&&&&&&&&&&&&&&&&&&&&&&\\
8 & 9 & 17 & 0 & 14 & 18 & 12 & 2 & 10 & 16 & 21 & 20 & 3 & 7 & 24 & 23 & 26 & 19 & 4 & 1 & 13 & 22 & 11 & 5 & 25 & 15 & 6\\
&&&&&&&&&&&&&&&&&&&&&&&&&&&\\
13 & 14 & 22 & 5 & 19 & 23 & 17 & 7 & 15 & 21 & 26 & 25 & 8 & 12 & 2 & 1 & 4 & 24 & 9 & 6 & 18 & 0 & 16 & 10 & 3 & 20 & 11\\
&&&&&&&&&&&&&&&&&&&&&&&&&&&\\
12 & 13 & 21 & 4 & 18 & 22 & 16 & 6 & 14 & 20 & 25 & 24 & 7 & 11 & 1 & 0 & 3 & 23 & 8 & 5 & 17 & 26 & 15 & 9 & 2 & 19 & 10\\
&&&&&&&&&&&&&&&&&&&&&&&&&&&\\
19 & 20 & 1 & 11 & 25 & 2 & 23 & 13 & 21 & 0 & 5 & 4 & 14 & 18 & 10 & 7 & 10 & 3 & 15 & 12 & 24 & 6 & 22 & 16 & 9 & 26 & 17\\
&&&&&&&&&&&&&&&&&&&&&&&&&&&\\
0 & 1 & 9 & 19 & 6 & 10 & 4 & 21 & 2 & 8 & 13 & 12 & 22 & 26 & 16 & 15 & 18 & 11 & 23 & 20 & 5 & 14 & 3 & 24 & 17 & 7 & 25\\
&&&&&&&&&&&&&&&&&&&&&&&&&&&\\
20 & 21 & 2 & 12 & 26 & 3 & 24 & 14 & 22 & 1 & 6 & 5 & 15 & 19 & 9 & 8 & 11 & 4 & 16 & 13 & 25 & 7 & 23 & 17 & 10 & 0 & 18\\
\end{smallmatrix}
\]
\vspace{.01in}
\begin{center}
\line(1,0){350}
\end{center}
\end{figure}

\subsection{An example of an encryption using E-DES}

We now describe how E-DES encrypts the $18$-nit message $m = 012012012012012012$ by using the $20$-nit key $k = 11012012122012012110$. 

{\flushleft{\bf Apply initial permutation $P$ to $m$:}}
\[
  m_1 = P(m) = 121020110102200221.
\]
The right half of $m_1$ is 
\[
  m_2=R(m_1) = 102200221.
\] 
{\flushleft{\bf Apply the expansion map $E$ to $m_2$:}}
\[
  m_3 = E(m_2) = 110220022002211.
\]
{\flushleft{\bf Apply the key compression map $CP_1$ to $k$:}}
\[
  k_1 = CP_1(k) = 120002201111211.
\]
{\flushleft{\bf In $\integers_3^{15}$ add $m_3$ and $k_1$:}}
\[
  m_4 = m_3\oplus k_1 = 110220022002211 \oplus 120002201111211 = 200222220110122.
\]                     
Note that the group operation $\oplus$ here is nit-wise addition modulo 3. 
{\flushleft{\bf Partition $m_4$ into three 5-nit blocks before processing to the S-boxes.}} Then the block $20022$ becomes an input in S-box 1, $22201$ an input in S-box 2 and the last 5-nit block $10122$ becomes an input in S-box 3.   
{\flushleft{\bf Determine the 3-nit output of each S-box:}}

{\flushleft{\bf S-box 1:}}
The first and last nits of this block form $22$, which corresponds to row number $8$, the ninth row. The middle three nits form $002$, which corresponds to column number $2$, which is the third column. The entry in row $8$, column $2$ of S-box 1 is the decimal number $20$. The base $3$ representation of $20$ is $202$, so the output for S-box 1 is $202$.
 
{\flushleft{\bf S-box 2:}}
The second S-box input is $22201$. The first and last nits of this block form $21$, which corresponds to row number $7$, which is the eight row of S-box 2. The middle three nits form $220$, which corresponds to column number $24$ which is the twenty-fifth column of S-box 2. The entry in row $7$, column $24$ of S-box 2 is the decimal number $5$. The base $3$ representation of $5$ is $012$, so the output for S-box $2$ is $012$. 

{\flushleft{\bf S-box 3:}}
The third S-box input is $10122$. The first and last nits of this block form $12$, which corresponds to row number $5$, which is the sixth row of S-box 3. The middle three nits form $012$, which corresponds to column number $5$, which is the sixth column of S-box 3. The entry in row $5$, column $5$ of S-box 3 is the decimal number $19$. The base $3$ representation of $22$ is $211$, so the output for S-box $3$ is $211$. 

{\flushleft{\bf Determine the 12-nit output from the S-boxes:}}

Concatenating these three S-box outputs in order gives the combined output
\[
  m_5 = 202012211.
\]
{\flushleft{\bf In the group $\integers_3$, add $m_5$ to the left half of $m_1$:}}
The left half of $m_1$ is $m_6 = 121020110$. This step gives
\[                               
  m_7 = m_6 \oplus m_5 = 020002021
\]
{\flushleft{\bf Combine $m_2$ and $m_7$ by left-right swap:}}
\[
  e_1 = \Theta(m_7,m_2) = 102200221020002021 
\]
This completes the first of the two Feistel rounds of E-DES.
The right half of $e_1$ is $e_2=R(e_1) = 020002021$.
{\flushleft{\bf Apply the expansion map $E$ to $e_2$:}}
\[
  e_3 = E(e_2) = 102000200020210.
\]
{\flushleft{\bf Apply the key compression map $CP_2$ to $k$:}}
\[
  k_2 = S_2(k) = 011010121202202.
\]
{\flushleft{\bf In $\integers_3^{15}$ add $e_3$ and $k_2$:}}
\[
  e_4 = e_3\oplus k_2 = 102000200020210 \oplus 011010121202202 = 110010021222112.
\]
{\flushleft{\bf Partition $e_4$ into three 5-nit blocks.}} Then the block $11001$ becomes an input in S-box 1, $00212$ an input in S-box 2 and the last 5-nit block $22112$ becomes an input in S-box 3.
{\flushleft{\bf Determine the 3-nit output of each S-box:}}

{\flushleft{\bf S-box 1:}}
The first and last nits of this block form $11$, which corresponds to row number $4$, the fifth row. The middle three nits form $100$, which corresponds to column number $9$, which is the tenth column. The entry in row $4$, column $9$ of S-box 1 is the decimal number $2$. The base $3$ representation of $2$ is $002$, so the output for S-box 1 is $002$.
 
{\flushleft{\bf S-box 2:}}
The second S-box input is $00212$. The first and last nits of this block form $02$, which corresponds to row number $2$, which is the third row of S-box 2. The middle three nits form $021$, which corresponds to column number $7$ which is the eighth column of S-box 2. The entry in row $2$, column $7$ of S-box 2 is the decimal number $8$. The base $3$ representation of $8$ is $022$, so the output for S-box 2 is $022$. 

{\flushleft{\bf S-box 3:}}
The third S-box input is $22112$. The first and last nits of this block form $22$, which corresponds to row number $8$, which is the ninth row of S-box 3. The middle three nits form $211$, which corresponds to column number $22$, which is the twenty-third column of S-box 3. The entry in row $8$, column $22$ of S-box 3 is the decimal number $23$. The base $3$ representation of $23$ is $212$, so the output for S-box $3$ is $212$. 

{\flushleft{\bf Determine the 12-nit output from the S-boxes:}}

Concatenating these three S-box outputs in order gives the combined output
\[
  e_5 = 002022212.
\]
{\flushleft{\bf In the group $\integers_3$, add $e_5$ to the left half of $e_1$:}}
The left half of $e_1$ is $e_6 = 102200221$. This step gives
\[                                
  e_7 = e_6 \oplus e_5 = 101222100.
\]
{\flushleft{\bf Concatenate $e_7$ and $e_2$ to form $e_7e_2$:}}

\[
  e_8 = e_7e_2 = 101222100020002021.
\]
{\flushleft{\bf Apply the final permutation, $P^{-1}$ to $e_8$:}}
\[
  c = P^{-1}(e_8) = 210212002210210000.
\]

Now $c=210212002210210000$ is the ciphertext output when E-DES is applied to the input plaintext $m = 012012012012012012$, using the key $k = 11012012122012012110$.

\section{Cycling closure experiments on simplified Feistel networks over certain finite groups}

We give a general overview of the cycling closure test. The cycling closure test can be used to address various questions about the algebraic structure of any finite cryptosystem. The test was used in \cite {CW} and in \cite {CO} to give a conclusive proof that the set of permutations in the classical DES is not closed and to determine the lower bound on the size of the subgroup generated by the DES permutations. We performed this test on several simplified versions of DES over certain finite groups. 

Assume that the subset $\TT_{\Pi}$ of the group ${\sf Sym}_{G^{2t}}$ is closed under composition. Then for any key $k\in \KK$ the order of the encryption permutation $T_{k}$ divides the order of $\TT_{\Pi}$. Recall that $\vert \TT_{\Pi}\vert\le \vert\KK\vert$. Thus, in particular, the order of $T_{k}$ is no larger than $\vert\KK\vert$. By the Orbit-Stabilizer theorem it follows that for any message $x\in\mathcal M$ we have $\vert {\sf orb}_{\langle T_k\rangle}(x)\vert$ divides the order of the cyclic group $\langle T_k\rangle$, and thus 
\[
  \vert {\sf orb}_{\langle T_k\rangle}(m)\vert \le \vert\KK\vert.
\]
More generally, for messages $x_1,\, x_2,\, \cdots,\, x_m$ we have
\[
  {\sf lcm}\{\vert {\sf orb}_{\langle T_k\rangle}(x_i)\vert:\, i\le m\} \le \vert\KK\vert.
\]
If it happens that for a key $k$ and messages $x_1,\,\cdots,\,x_m$ 
\[
{\sf lcm}\{\vert {\sf orb}_{\langle T_k\rangle}(x_i)\vert:\, i\le m\} > \vert\KK\vert
\]
then we have that $\TT_{\Pi}$ is not closed and therefore not a group. 

Thus to test a cryptosystem for some algebraic weaknesses such as closure and purity one has to examine the orbits of subsets of encryption transformations on particular messages. The method is to compute the orbits of single encryption  and to apply the cycling closure test to subsets of two or more encryption transformations. 
The cycling closure test picks an initial message $m$ at random and then takes a pseudorandom walk in $\TT_{\Pi}$, beginning at $m$. For each step of the pseudorandom walk, the previous ciphertext is encrypted under a key chosen by a pseudorandom function of the previous ciphertext. The walk continues until a cycle is detected. By the Birthday Paradox the walk is expected to cycle after approximately $\mid \TT_{\Pi}\mid^{1/2}$ steps.

Recall the definition of purity of a cryptosystem. To determine the purity of the cryptosystem using the cycling closure test we first need show the following statement.

\begin{lemma}
A cryptosystem $\Pi$ is pure if and only if for some $T\in\TT_{\Pi}$ the set $T^{-1}\TT_{\Pi} = \{T^{-1}\circ T_k:k\in\KK\}$ is closed under functional composition.
\end{lemma}
\begin{pf}
Let $T\in\TT_{\Pi}$ be an encryption permutation for which the set $T^{-1}\TT_{\Pi}$ is closed under functional composition. Then, for each $T_i,T_j\in \TT_{\Pi}$ and each $T\in\TT_{\Pi}$ there exist $T_k\in \TT_{\Pi}$ such that $(T^{-1}\circ T_i)\circ (T^{-1}\circ T_j)=T^{-1}\circ T_k$. But, then 
\[
T\circ(T^{-1}\circ T_i)\circ (T^{-1}\circ T_j)=T\circ T^{-1}\circ T_k 
\]
or $T_i\circ T^{-1}\circ T_j= T_k$ which means that $\Pi$ is pure.

Now, let assume that $\Pi$ is pure. Then for each  $T_i,T_j, T_k\in \TT_{\Pi}$ there exists $T_l\in \TT_{\Pi}$ such that  $T_i\circ T_j^{-1}\circ T_k=T_l$. But, then for each $T\in \TT_{\Pi}$, 
\[
(T^{-1}\circ T_i\circ T_j^{-1})\circ (T^{-1}\circ T_k)=T^{-1}\circ T_l
\]
This implies that for some $T\in \TT_{\Pi}$ the set $T^{-1}\TT_{\Pi}$ is closed. $\epf$
\end{pf}

\subsection{\bf Orbit Test.} Given any key $k$ and any message $m$, compute $x_i=T_k^{i}(m)$, $i=1,2,\cdots$ for a specified numbers of steps or until a cycle is detected.

\subsection{\bf Purity Test.} Pick any encryption $T\in \TT_{\Pi}$ and apply the cycling closure test to the set $T^{-1}\circ\TT_{\Pi}$.
\subsection{\bf Small Subgroup Test.} Pick any keys $k_1,k_2,\cdots k_s$ and any message $m$ apply the cycling closure test to the set $\{T_{k_1},T_{k_2},\cdots, T_{k_s}\}$  $T\in \TT_{\Pi}$ to obtain a statistical lowerbound of the group $\langle T_{k_1},T_{k_2},\cdots, T_{k_s}\rangle$. \\ 

We developed a software to implement these tests on a 2-round simplified version of DES with two randomly chosen S-boxes where the group operation is addition modulo $n$ for $n\in \{2,3,5,7,11\}$. Tables 2-4 give a through description of our cycling experiments. The tests give a conclusive proof to the following theorems.

\begin{theorem}\label{coppersmith}
For a fixed pair of S-boxes and for $n\in \{2, 3, 5, 7, 11\}$ the set of encryption permutations of a 2-round simplified version of DES over $\integers_n$ is not closed under functional composition.
\end{theorem}
  
\begin{theorem}
The simplified cipher E-DES is not pure.
\end{theorem}

\begin{theorem}
The group generated by the \emph {E-DES} encryptions is larger than $\vert\sym_{49}\vert$.  
\end{theorem}

\section{Some group theoretic properties of\\ Feistel based block ciphers}

In \cite {LR} it was shown that the Luby-Rackoff constructions with 4 rounds based on the XOR operation are secure against adaptive chosen plaintext and ciphertext attacks. These results are based on the rather strong hypothesis that the round functions are random. In \cite {P} it was shown that when the round functions are random permutations, a 4-round Feistel network remains secure as long as the number of queries $m$ is very small compared with $2^{t/2}$ (i.e. $m\ll 2^{t/2}$) where $2t$ is the block size. One way to improve the security of these type of Feistel networks is to use multiple encryptions. 

We show that the $n$-round Feistel networks, $n\leq 6$ with one-to-one round functions `` do not form a group" i.e. the set of such Feistel permutations do not form a group under functional composition. This implies that multiple encryptions can improve the security of the 4-round Feistel networks even when the round functions are random permutations. 

\begin{theorem}
Let $G$ be a finite group with $\vert G\vert >1$ and let $t$ be a positive integer. Let $X\subset {\sf F}_t(G)$ be a set of functions $f:G^t\rightarrow G^t$ that does not include the identity element. If each element of $X$ 
is one-to-one then for $n\leq  6$ the set of permutations of the form $\Psi^{2t}(f_1,f_2,\cdots, f_n)$ where each $f_i,\, i\le n$ is from $X$ is not a subgroup of $\sym_{\mid G^{2t}\mid}$. 
\end{theorem}
\begin{pf}
We prove the statement for $n=6$. The proof for $n<6$ is similar. Consider the set of all permutation of the form $\Psi^{2t}(f_1,f_2,\cdots, f_6)$. Assume that this set is a group under composition. Then some $\epsilon$ must be the identity element of $\sym_{\mid G^{2t}\mid}$. This means that for some fixed $f\in X$, for all $(x,y)\in G^t\times G^t$, we have  $\Psi^{2t}(f_1,f_2,\cdots, f_6)(x,y)=(x,y)$. But then for all $x, y \in G^t$ we have 
\pagebreak
\begin{equation}\label{fx}
f_1(y)+f_3(y+f_2(x+f_1(y)) +f_5(y+f_2(x+f_1(y))+f_4(x+f_1(y)+
\end{equation}
\[
+ f_3(y+f_2(f+f_1(y))))=0
\]
and 
\begin{equation}\label{fy}
f_2(x+f_1(y))+f_4(x+f_1(y)+f_3(y+f_2(x+f_1(y)))+f_6(y)=0
\end{equation}
Put $x=-f_1(y)$ in \ref{fx} and \ref{fy} to get
\begin{equation}\label{fx1}
f_1(y)+f_3(y+f_2(0))+f_5(y+f_2(0)+f_4(f_3(y+f_2(0)))=0
\end{equation}
and
\begin{equation}\label{fy1}
f_2(0)+f_4(f_3(y+f_2(0)))+f_6(y)=0
\end{equation}
After substituting \ref{fy1} into \ref{fx1} and then setting $y=0$ we get
\[
f_1(0)+f_3\circ f_2(x+f_1(0))+f_5\circ (-f_6)(0)=0
\]
Thus, 
$f_3\circ f_2$ is a constant function. Since $f_2\in {\sf F}_t(G)$ is one-to-one function one concludes that $f_3$ must be a constant. But this contradicts the fact that $f_3$ is a one-to-one function\footnote{Observe that the theorem can be strengthened by requiring only that $f_2$ is a one-to-one member of ${\sf F}_t(G)$, and $f_3$ is a non-constant element of ${\sf F}_t(G)$.}.
$\epf$
\end{pf}
\vspace{0.1in}

Generalizing classical DES, we now define for the finite group $G$ and positive integers $t$ and $n$, $n$-round DES over $G$, denoted $\mbox{GDES}^n_{2t}$: For given functions $f_1,\,\cdots,\, f_n$ the corresponding $n$-round $\mbox{GDES}^n_{2t}$ permutation $T_n$ over the finite group $(G,\oplus)$ is the composition 
\[
   T_n = P^{-1}\circ\theta\circ\sigma_{f_n}\circ\sigma_{f_{n-1}}\circ\cdots\circ\sigma_{f_1}\circ P
\] 
of permutations, where
$\theta$ is the ``swap" function 
           \begin{center}
               $\theta(x,y)=(y,x)$
            \end{center}
from $G^{2t}$ to $G^{2t}$ and $P$ is a member of $\sym_{\mid G^{2t}\mid}$ and $P^{-1}$ its inverse. These are the initial and final permutations used in $\mbox {GDES}^n_{2t}$. In this notation classical DES is ${\integers}_2$DES$^{16}_{64}$.

Note that $\theta = \theta^{-1}$. In the case when the underlying group $G$ in $\mbox{GDES}^n_{2t}$ is a group of characteristic $2$, $\sigma_f = \sigma_f^{-1}$. For arbitrary finite groups, 
\[
\sigma_f^{-1}(x,y) = (x \ominus f(y), y),
\]
where $\ominus$ is the inverse group operation $x \ominus y = x \oplus (-y)$. Hence, in general, the decryption process applies the key schedule in the reverse order, with $\ominus$ used instead of $\oplus$; that is,
\[
T_n^{-1} = P^{-1} \circ \sigma_{f_1}^{-1}\circ \cdots \circ  \sigma_{f_n}^{-1}\circ\theta \circ P.
\]

Also, note that if the round functions $f_1,\,\cdots,\, f_n$ are randomly and independently chosen functions then
\[
T_n=P^{-1}\circ\theta\circ \Psi^n_{2t}\circ P
\] 

For a fixed set $X\subseteq {\sf F}_t(G)$ of functions, $\mbox {GDES}^n_{2t}(X)$ denotes $\mbox{GDES}^n_{2t}$ where the functions used in constructing the $\mbox{GDES}^n_{2t}$ permutations are restricted to come from $X$. We consider the effect of structure of subsets $X$ of $(({\sf F}_t(G),\, \odot),\,\circ)$ on the corresponding set of $\mbox{GDES}^n_{2t}(X)$ permutations. 
Note that the initial and final permutation, $P$ and $P^{-1}$ do not affect the algebraic structure of $\mbox{GDES}^n_{2t}$ and therefore they can be omitted when examining its algebraic properties.
 
\begin{theorem}\label{smallroundsth} Let $G$ be a finite group.
Let $X\subseteq {\sf F}_t(G)$ be a set of functions that does not include the identity element of $({\sf F}_t(G),\odot)$.  
Then 
\begin{enumerate}
  \item{The set of \emph{GDES}$^1_{2t}(X)$ permutations is not a subgroup of $\sym_{\mid G^{2t}\mid}$.}
  \item{The set of \emph{GDES}$^2_{2t}(X)$ permutations is not a subgroup of $\sym_{\mid G^{2t}\mid}$.}
 \end{enumerate}
\end{theorem}
\begin{pf}

Consider $\mbox{GDES}^1_{2t}(X)$: Assume that the set $\{\epsilon(x,y)=(x\oplus f(y),y): f\in X\}$, is a group under composition. Then some $\epsilon$ must be the identity element of $\sym_{\mid G^{2t}\mid}$. This means that for some fixed $f\in X$, for all $(x,y)\in G^t\times G^t$, we have  $(x\oplus f(y),y)=(x,y)$. But then for all $y$ in $G^t$, $f(y)= e$, the identity element of $G^{t}$, whence $f$ is the identity element of $({\sf F}_t(G),\odot)$. This contradicts the fact that $X$ does not include the identity element of $({\sf F}_t(G),\odot)$.

Next, consider $\mbox{GDES}^2_{2t}(X)$: Assume that the set of all $\mbox{GDES}^2_{2t}(X)$ permutations is a group under composition. Then there are functions $f_1,\, f_2\in X$ such that $\epsilon=\theta\circ\sigma_{f_2}\circ\sigma_{f_1}$ is the identity element of $\sym_{\mid G^{2t}\mid}$. Then for each  $(x,y)\in G^t\times G^t$, 
\begin{center}
\begin{tabular}{lcl}
$x$ & = & $x\oplus f_1(y)$ and \\
$y$ & = & $y\oplus f_2(x\oplus f_1(y)).$ 
\end{tabular}
\end{center}

But then for all $x$ and $y$ in $G^{t}$, $f_{1}(y)=f_{2}(x)= e $, the identity element of $G^t$, whence $f_1$ and $f_2$ are the identity element of $({\sf F}_t(G),\odot)$. This contradicts the fact that $X$ does not include the identity element of $({\sf F}_t(G),\odot)$.
$\epf$
\end{pf}\\
Clearly this deductive argument improves Theorem \ref{coppersmith} which was proven using a computational method.
\begin{corollary}\label{variousZn}
For any pair of S-boxes and for $n\in \{2, 3, 5, 7, 11\}$ the set of encryption permutations of a 2-round simplified version of DES over $\integers_n$ is not closed under functional composition.
\end{corollary}

Using similar techniques we have shown that the set of $\mbox{GDES}^n_{2t}(X)$ permutations for $n<6$ do not constitute a group under functional composition. 
\begin{theorem}
Let $G$ be a finite group and let $t$ be a positive integer. If $X\subseteq {\sf F}_t(G)$ is  the set of all one-to-one functions $X\subseteq {\sf F}_t(G)$ then the set of $\mbox{GDES}^n_{2t}(X)$ permutations for $n<6$ is not a subgroup of $\sym_{\mid G^{2t}\mid}$ 
\end{theorem}

The following corollary is an immediate consequence of the previous two theorems.

\begin{corollary}
For each of \emph{B-DES}, \emph{S-DES} and \emph{E-DES}, the corresponding set of encryption permutations does not constitute a group under functional composition.
\end{corollary}

\begin{definition}
The {\it characteristic} of the group $(G,\oplus)$ is the least positive integer $n$ such that for each $x\in G$ we have $nx=e$.  
\end{definition}

Assuming that the group has characteristic 2 and that the set of $\mbox{GDES}^{n-2}_{2t}$ permutations does not contain the identity element of $\sym_{\mid G^{2t}\mid}$ we can extend the previous result to $n\geq 6$ (Theorem \ref {gdes}). This partially answers the question from \cite{Kaliski} whether the set of encryption permutations in classical DES contains the identity element. Note that showing that the set of $\mbox{GDES}^{n}_{2t}$ permutations does not contain the identity element implies that it is also not closed under functional composition.

\begin{theorem}\label{gdes}
Let $G$ be a finite group of characteristic 2. If for each instance \emph{GDES}$^{n-2}_{2t}$ the set of \emph{GDES}$^{n-2}_{2t}$  encryption permutations does not contain the identity element of $\sym_{\mid G^{2t}\mid}$, then  for each instance of \emph{GDES}$^n_{2t}$ the subset of \emph{GDES}$^n_{2t}$ encryption permutations for which $f_1 = f_n$ is not a subgroup of $\sym_{\mid G^{2t}\mid}$.
\end{theorem}
\begin{pf} 

Recall that\footnote{Note that the initial and the final permutation don't affect the algebraic structure of $\mbox{GDES}^n_{2t}$} $\mbox{GDES}^n_{2t}$ permutations are of the form 
\[
  \epsilon_n=\theta\circ\sigma_{f_n}\circ\sigma_{f_{n-1}}\circ\cdots\circ\sigma_{f_2}\circ\sigma_{f_1}.
\]
Let $(x_i,y_i)$ denote $\epsilon_{i}(x,y)$ for $(x,y)\in G^t\times G^t$ and $1\leq i\leq n$. Assume that the subset of $\mbox{GDES}^n_{2t}$ permutations for which $f_1 = f_n$ contain the identity element of $\sym_{\mid G^{2t}\mid}$.

Fix an instance of $\mbox{GDES}^n_{2t}$ for which there is a key giving rise to the sequence $(f_1,\cdots,f_n)$ of round functions such that $f_1=f_n$, and the corresponding $\mbox{GDES}^{n}_{2t}$ encryption permutation is the identity function. Also fix a key $k$ giving rise to such a sequence $(f_1,\cdots,f_n)$ of round functions. Thus, for all $x,\, y\in G^t\times G^t$ we have $\epsilon_n(x,y) = (x,y)$.

Note that from this instance of $\mbox{GDES}^n_{2t}$ we can define an instance of $\mbox{GDES}^{n-1}_{2t}$ so that the key schedule of the latter is related as follows to the key schedule of the former: For the given key $k$, if $(k_1,\cdots,k_n)$ are the $n$ round keys of $\mbox{GDES}^n_{2t}$, then $(k^{\prime}_1,\cdots,k^{\prime}_{n-1})$ are the corresponding round keys for $\mbox{GDES}^{n-1}_{2t}$, where $k^{\prime}_i = k_{i+1}$ for $1\le i < n$. 

Then we can write $\epsilon_n=\epsilon^{\prime}_{n-1}\circ\sigma_{f_1}$ where $\epsilon^{\prime}_{n-1}=\theta\circ\sigma_{f_n}\circ\sigma_{f_{n-1}}\circ\sigma_{f_3}\circ...\circ\sigma_{f_2}$ is the corresponding $\mbox{GDES}^{n-1}_{2t}$ encryption permutation arising from the key $k$. For convenience, write $(x_{n-1},y_{n-1})$ for $\epsilon^{\prime}_{n-1}(x,y)$. 

Since we assumed that $\epsilon_n$ is the identity function we find that  $\epsilon^{\prime}_{n-1}= \sigma_{f_1}^{-1}$, \emph{i.e.}, for all $x$ and $y$ in $G^t$,
\begin{equation}\label{nminoneeq}
\epsilon^{\prime}_{n-1}(x,y)=(y\ominus f_{1}(x),x)
\end{equation}

On the other hand we can write $\epsilon^{\prime}_{n-1}(x,y) = \theta\circ \sigma_{f_n} \circ \theta \circ \epsilon^{\prime}_{n-2}(x,y)$, where as before $\epsilon^{\prime}_{n-2}$ is a corresponding $\mbox{GDES}^{n-2}_{2t}$ encryption permutation obtained from the key $k$ using the ideas above to define an instance of $\mbox{GDES}^{n-2}_{2t}$ from the instance of $\mbox{GDES}^{n-1}_{2t}$. Note that for all $x$ and $y$ in $G^t$,
\begin{equation}\label{nmintwoeq}
  \theta\circ\sigma_{f_n}\circ\theta(x,y) = (y\oplus f_{n}(x),x). 
\end{equation}

From equations (\ref{nminoneeq}) and (\ref{nmintwoeq}) we see that for all $x$ and $y$ in $G^t$ we have, setting $\epsilon^{\prime}_{n-2}(x,y) = (x_{n-2},y_{n-2})$, that
\[
  (y\ominus f_1(x),x) = (y_{n-2}\oplus f_n(x_{n-2}),x_{n-2}).
\]
Thus, for all $x$ and $y$ in $G^t$ we have that $x_{n-2} = x$, and $y\ominus f_1(x) = y_{n-2}\oplus f_n(x)$. Since we assumed that $f_1 = f_n$ and that the group has characteristic $2$, we find that for all $y\in G^t$, $y_{n-2}=y$.
But then the encryption permutations of this instance of $\mbox{GDES}^{n-2}_{2t}$ contain the identity element of $\sym_{\mid G^{2t}\mid}$. This establishes the contrapositive of the theorem. $\epf$
\end{pf}

\vspace{0.1in}

In \cite {BBKW} we show when the group generated by the $n$-round $\mbox{GDES}^n_{2t}$ permutations, where the underlying Feistel network  based on the binary operation of $G$ contains odd permutations.

\begin{theorem}{\bf \cite {BBKW}}
Let $G$ be a  finite group and let $n$ be positive integer. If $\mid G\mid ^t \equiv 2,3\; mod \;4$ and $t$ is odd then the group generated by the set of \emph{GDES}$^n_{2t}$ encryption permutations is not a subgroup of $\alt_{\mid G^{2t}\mid}$.
\end{theorem}

\begin{corollary}
The subgroup of $\sym_{3^{18}}$ that is generated by the set of \emph{E-DES} permutations is a group that contains  as many even permutations as odd permutations. 
\end{corollary}

Based on several asymptotic results that appear in the literature ( \cite{B}, \cite{D}, or \cite {MT}) about generating the alternating group or the symmetric group we make the following conjecture. 

\begin{conjecture}
The group generated by the set of \emph{E-DES} permutations is $\sym_{3^{18}}$.
\end{conjecture}

\begin{theorem}\label{sboxexpansion}
Let $G$ be a subgroup of the finite group $H$. If there is an S-box set for $G$ such that the corresponding set of \emph{GDES}$^2_{2t}$ encryption permutations do not constitute a group, then there is an S-box set for $H$ such that the  corresponding set of \emph{HDES}$^2_{2t}$ encryption permutations do not constitute a group.
\end{theorem}

\begin{pf}
Let S-boxes $S_1,\cdots,S_n$ for $\mbox{GDES}^2_{2t}$ be given such that the set of encryption permutations of $\mbox{GDES}^2_{2t}$  defined from these S-boxes does not constitute a group under composition. We may assume that $t = j\hspace{-0.03in}\cdot\hspace{-0.03in} n$, while round key lengths are $(i+j)\hspace{-0.03in}\cdot \hspace{-0.03in}n$. Thus each of these S-boxes has $\vert G^i\vert$ rows and $\vert G^j\vert$ columns. Each of the row entries is an element of $G^j$. We may assume that each of the $\vert G\vert^i$ rows is indexed by an element of $G^i$ and each of the $\vert G\vert^j$ columns is indexed by elements of the set $G^j$.

Now expand each S-box $S_p$ to an S-box $R_p$ as follows: Adjoin $\vert H^i\vert - \vert G^i\vert$ additional rows, and $\vert H^j\vert - \vert G^j\vert$ additional columns to obtain an $\vert H^i\vert$-by-$\vert H^j\vert$ array such that the top left $\vert G^i\vert$-by-$\vert G^j\vert$ corner is the given S-box $S_m$, and the remaining entries are all from $H^j\setminus G^j$. The additional rows are now indexed by elements of $H^i\setminus G^i$, and the additional columns are indexed by elements of $H^j\setminus G^j$.

{\flushleft{\bf Claim: }} The set of $\mbox{HDES}^2_{2t}$ permutations arising from the S-boxes $R_1,\cdots,R_n$ does not constitute a group under functional composition.

For assume the contrary. Let $k_1$ and $k_2$ be $\mbox{GDES}^2_{2t}$ keys such that $T_{k_2}T_{k_1}$ is not of the form $T_{k_3}$for some $\mbox{GDES}^2_{2t}$ key $k_3$. Observe that the set of $\mbox{GDES}^2_{2t}$ keys is  subset of the set of $\mbox{HDES}^2_{2t}$ keys, and that the set of $\mbox{GDES}^2_{2t}$ plaintexts is a subset of the set of $\mbox{HDES}^2_{2t}$ plaintexts. 

By our hypothesis there is an $\mbox{HDES}^2_{2t}$ key $k$, fixed from now on, such that 
\[
  T_{k_2}T_{k_1} = T_k.
\]
By the choice of the keys $k_1$ and $k_2$ it follows that $k$ has an entry from $H\setminus G$. It also follows that there is a plaintext $m\in G^{2t}$ such that during encryption of $m$ using $T_k$, in some of the two Feistel rounds an entry of $k$ which is in $H\setminus G$ is used, for otherwise we may modify $k$ so that all entries are from $G$ and still have $T_{k_2}T_{k_1}=T_k$, contradicting the choice of $k_1$ and $k_2$.

{\flushleft{\underline{Encryption of $m$ by $T_{k_2}T_{k_1}$:}}}
Let $m_{\ell}$ denote the left half of $m$ and let $m_r$ denote the right half of $m$. Consider $E(m_r)$, the expansion of $m_r$, used in specifications of the encryption algorithm $\mbox{GDES}^2_{2t}$ in use. 

In round $1$ of $T_{k_1}$ we see that the first step is $m^1_1 = E(m_r) \oplus k^1_1$ where $k^1_1$ is the first round subkey. Now all $(i+j)\cdot n$ entries of $m^1_1$ are elements of $G$, and thus point to entries of the S-boxes $S_1,\cdots,S_n$, so that the output of this round is another element, $e^1_1$, of $G^{2t}$. Completing the second round of $T_{k_1}$ produces an element $e^1_2$ of $G^{2t}$, which is input of $T_{k_2}$. By similar considerations the result of applying $T_{k_2}$ to $e^1_2$ is the element $T_{k_2}T_{k_1}(m)$ of $G^{2t}$.

{\flushleft{\underline{Encryption of $m$ by $T_{k}$:}}}
Now consider the destiny of $m$ under the encryption permutation $T_k$. By our earlier remark about the use of entries of $k$ during the encryption process, for some round $j\le 2$, the round key $k^j$ has an entry in $H\setminus G$ which is used in that round. Let $j$ denotes the first round when this occurs. 
{\flushleft{\bf Case 1: j=1.}} 
The input to round $1$ is $m\in G^{2t}$. With $E(m_r)$ the expansion of the right half to an $(i+j)\hspace{-0.03in}\cdot\hspace{-0.03in} n$ sequence from $G$, $m^1 = E(m_r) \oplus k^1$ has an entry from $H\setminus G$ as it is a group element of $H$ obtained from the group operation of an element of the subgroup $G$ with an element not in $G$.  But then $m^1$ points, for some S-box $R_{p}$ either at a row beyond the $\vert G^i\vert$-th, or at a column beyond the $\vert G^j\vert$-th. By the construction of the expanded S-boxes over $H$, the return from this S-box consists of $j$-nits from $H^j\setminus G^j$, and thus the output from this step is $e_1$, which has $n\cdot j$ nits, and one of the $n$ blocks of consecutive nits has an entry from $H\setminus G$.
The output from this round then is of the form a list of length $2\hspace{-0.03in}\cdot \hspace{-0.03in}t$ with first $t$ entries $m_r$, the right half of the original message, a member of $G^t$, and $m_{\ell}\oplus e_1$ which is a member of $H^t\setminus G^t$. After the second round of the encryption using the key $k$ we find that the right half of the output ciphertext still is $m_{\ell}\oplus e_1$, and so the result of the encryption is not a member of $G^{2t}$, contradicting that for all $m$, $T_{k_2}T_{k_1}(m) = T_k(m)$.

{\flushleft{\bf Case 2: j=2.}}
The output of the first round of $T_k$ is an element of $G^{2t}$. Thus the input to round $2$ is an $M\in G^{2t}$. Let $M_{\ell}$denote the left $t$ nits of $M$, and $M_r$ the right $t$ nits. Let $E(M_r)$ be the expansion of the right half to an $(i+j)\hspace{-0.03in}\cdot\hspace{-0.03in} n$ sequence from $G$. Then $M^1 = E(M_r) \oplus k^2$ has an entry from $H\setminus G$ as it is a group element of $H$ obtained from the group operation of an element of the subgroup $G$ with an element not in $G$. 

But then $M^1$ points, for some S-box $R_{p}$ either at a row beyond the $\vert G^i\vert$-th, or at a column beyond the $\vert G^j\vert$-th. By the construction of the expanded S-boxes over $H$, the return from this S-box consists of $j$-nits including some from $H\setminus G$. Thus the return from the S-boxes is an element of $H^t\setminus G^t$, say $e$. But then $M_{\ell} \oplus e$ still is a member of $H^t\setminus G^t$ as $G$ is a subgroup of $H$. Since this is the right half of the output from the second round of the encryption using the key $k$, it follows that $T_k(m)\in H^{2t}\setminus G^{2t}$, and thus $T_{k_2}T_{k_1}(m)\neq T_k(m)$, a contradiction.
 
It follows that set of encryption permutations of $\mbox{HDES}^2_{2t}$ with the S-boxes $R_1,\cdots,R_n$ does not constitute a group.
$\epf$
\end{pf}

The following corollary follows directly from the previous theorem.

\begin{corollary}
Let $G$ be an arbitrary finite group. If there is an S-box set for \emph{GDES}$^2_{2t}$ encryption permutations do not constitute a group, then for any finite group $H$ there is an S-box set such that the set of \emph{UDES}$^2_{2t}$ encryption permutations where $U=G\times H$ do not constitute a group.
\end{corollary}

Note that the last $\vert H\vert^i - \vert G\vert^i$ rows of the S-boxes $R_1$, $\cdots$, $R_n$ appearing in the proof of Theorem \ref{sboxexpansion} do not contain any members of $G^j$, and thus are not permutations of members of $H^j$. We conjecture that expansions of the original set $n$ of S-boxes can be found such that if the rows of the original S-boxes were permutations of $G^j$, then the rows of the expanded S-boxes will be permutations of $H^j$, and the theorem would hold.

\section {Conclusions and future work}
This paper is a study of DES-like ciphers over arbitrary finite groups. We introduced a new simplified version of DES based on the finite group ${\mathbb Z}_3$. We showed that its set of encryption permutations does not form a group under functional composition. Corollary \ref{variousZn} indicates that this result is not particular to the underlying group $\integers_3$. Theorem \ref{smallroundsth} shows that this conclusion holds for arbitrary finite groups. 
Our proof of this fact deviates from former computationally intensive methods by being a purely deductive proof. Before we discovered our deductive proof we also used the cycling closure test as was used in the case of DES for the same question. All examples of DES-like cryptosystems are based on commutative groups. Though Theorem \ref{smallroundsth} applies also to this case, we have not yet explored potential new complications that may arise when a version of DES is based on a non-commutative group.

The group generated by this simplified DES cipher is not the alternating group as opposed to the group generated by the one-round functions of the classical DES and other block ciphers that generate the alternating group. It would be interesting to determine which finite groups can be the group generated by the set of encryption permutations of a DES-like cipher. 

We showed that the $n$-round Feistel networks, $n\leq 6$ with one-to-one round functions ``do not form a group" i.e. the set of such Feistel permutations do not form a group under functional composition. It will be useful to know whether one can extend this result to hold for any number of rounds and based on any group operation in the underlying Feistel network. Our results motivate a need for re-examining the DES-like ciphers to determine the extent to which the old results hold when we consider arbitrary finite algebraic structures.

\vspace{0.05in}

\vskip 5truept

\begin{flushleft}
$^1$ \small{Department of Mathematics, Boise State University, Boise ID 83725}\\
$^2$ \small {Department of Mathematics, Loyola Marymount University, Los Angeles, CA 90045}\\
$^3$ \small {Department of Mathematics, University of Western Carolina, Cullowhee, NC 28723}
\end{flushleft}

\vspace{4in}
\appendix
\section{Cycling Experiments}

{\bf Orbit Test}: The closure experiment for ${\integers}_2$DES$^{2}_8$, ${\integers}_3$DES$^{2}_8$, ${\integers}_7$DES$^{2}_8$,  and ${\integers}_{11}$DES$^{2}_8$ with random S-boxes, initial and final permutation and expansion function that meet their standard architectural requirements of Feistel based cryptosystems.
\vspace{0.3in}
\begin{center}
\begin{tabular}{|c|lclllc||l|}
	\hline
\multicolumn{8}{|c|}
	{\rule[-3mm]{0mm}{8mm} {\bf Table 2}: Computing $\mid\{T_{k_j}^r(m_i): r=1,2,\cdots\}\mid $}\\ \hline
\text{$G$}
                 & {\small $m_i$} &              & &  {\small $k_j$}         &{\small $orb(m_i)$} &    &{\small $lcm(orb(m_i):i\in\naturals)$}           \\ \hline\hline 
$\integers_2$    & 130       &              & &   823                          &   31     &    &                   \\
                 & 212       &              & &   963                          &   37     &    & $1147 > 2^8$      \\ \hline
$\integers_3$    & 5838      &              & & 41155                          & 2526     &    &                   \\     
                 & 5580      &              & & 37372                          & 1739     &    & $4392714 > 3^8$    \\ \hline
$\integers_5$    & 177954    &              & & 2817216                        & 8350     &    &                   \\
                 & 240687    &              & & 6114274                        & 46728    &    & $195089400 > 5^8$        \\ \hline
$\integers_7$    & 4903806   &              & & 161418036                      & 1377440  &    &                   \\
                 & 4684968   &              & & 104555592                      & 3014559  &    & $4152374148960 > 7^8$    \\ \hline
$\integers_{11}$ & 40015435  &              & & 4533979344                     & 106572673&    &                   \\
                 & 110072914 &              & & 25730291171                    & 19064231 &    & $2031726056359463 > 11^8$\\ \hline
\end{tabular}
\end{center}
\vspace{0.3in}
\begin{center}
\begin{tabular}{|lclllllllllc||l|}
	\hline
\multicolumn{13}{|c|}
	{\rule[-3mm]{0mm}{8mm} {\bf Table 3}: Purity experiment for E-DES}            \\ \hline
      {\small $m$} &&&&&  {\small $e_K$}    &&& {\small $d_K$} &&&&  {\small $orb(m)$}                                                                                        \\ \hline\hline 
         67681038  &&&&&   22933471         &&&   1402043471   &&&&      12802413       \\
                 
         884024783 &&&&&   1402043471       &&&   9625730992   &&&&      254125864      \\     \hline

\end{tabular}
\end{center}

\vspace{0.3in}

\begin{center}
\begin{tabular}{|clclllllllllllc||llll|c||c||l|l|} \hline
        \multicolumn{19}{|c|}
	{\rule[-3mm]{0mm}{8mm} {\bf Table 4} : Small subgroup experiment for E-DES}             \\ \hline
             &&{\small $m_i$} &      &&&&   &     &  {\small $k_j$} &&            &&&&{\small $orb(m_i)$} &&& \\ \hline\hline 
             &&     0         &      &&&&   &     & 4533979344      &&            &&&& 134282729    &&&            \\ 
             &&     0         &      &&&&   &     & 1402043471      &&            &&&& 216589023   &&&            \\
             &&     6021247   &      &&&&   &     & 1402043471      &&            &&&& 201375970   &&&            \\
             &&     746014783 &      &&&&   &     &   0             &&            &&&& 62909599    &&&            \\
             &&     580027391 &      &&&&   &     &   0             &&            &&&& 201375970   &&&            \\
             &&     442017391 &      &&&&   &     & 1402043471      &&            &&&& 134282729   &&&            \\ 
             &&     746014783 &      &&&&   &     & 1402043471      &&            &&&& 18939453    &&&            \\
             &&     872037801 &      &&&&   &     & 1402043471      &&            &&&& 68600442    &&&            \\
             &&     67681038  &      &&&&   &     & 1402043471      &&            &&&& 134282729   &&&            \\\hline 

\end{tabular}
\end{center}

\end{document}